\documentclass[12pt]{article}
\usepackage{amssymb}

\def\triangulo{\hbox{.\kern2pt.\kern-5pt\raise 4pt\hbox{.}}\kern 5pt}
\font\gorditas = msbm8
\def\bbb#1{\hbox {{\gordas #1}}}
\def\errita{\hbox{\gorditas R}}

\font\gordas = msbm10 at 12pt
\def\bbb#1{\hbox {{\gordas #1}}}
\def\a{{\bbb A}}
\def\erre{{\bbb R}}
\def\ce{{\bbb C}}
\def\e{{\bbb E}}
\def\o{{\bbb O}}
\def\ze{{\bbb Z}}
\def\ache{{\bbb H}}

\begin{document}
\begin{center}{\bf Hopf construction map in higher dimensions}\\[1cm]
Guillermo Moreno
\end{center}
\vglue1cm
\noindent
{\bf Abstract:} In this paper, we study the zero set of the Hopf construction map
$F_n:\a_n\times \a_n\rightarrow\a_n\times \a_0$
given by $F_n(x,y)=(2xy,||y||^2-||x||^2)$ for $n\geq 4$, where $\a_n$ is the Cayley-Dickson
algebra of dimension $2^n$ on $\erre$.

\footnote{\kern-.7cm Keywords and phrases: Cayley-Dickson algebras, alternative 
algebras, zero divisors, flexible algebra,normed algebra.}

\vglue1cm
\noindent
{\bf Introduction:} Let  $f_1:S^3\rightarrow S^2, f_2:S^7\rightarrow S^4$ and $f_3:S^{15}\rightarrow S^8$ be the classical {\bf Hopf maps}; these can be defined using {\bf the Hopf construction}: Let $\a_1=\ce, \a_2=\ache$ and $\a_3=\o$ be the complex, quaternion and octonion numbers respectively and $F_n:\a_n\times \a_n\rightarrow\a_n\times\erre$ are given by
$$
F_n(x,y)=(2xy,||y||^2-||x||^2)$$
for $n=1,2,3$.

Now write $S^{2^{n+1}-1}=\{(x,y)\in \a_n\times\a_n:||x||^2+||y||^2=1\}$.
 By definition, 
$$F_n|S^{2^{n+1}-1}=f_n$$ are the Hopf maps. Since $\a_n$ 
is a normed real algebra of dimension $2^n$,  for $n=1,2,3$ we have that
\begin{eqnarray*}
||(2xy,||y||^2-||x||^2)||^2&=&4||xy||^2+(||y||^2-||x||^2)^2\\
&=&4||x||^2||y||^2+||y||^4+||x||^4-2||x||^2||y||^2\\
&=&(||x||^2+||y||^2)^2,
\end{eqnarray*}
so if $||x||^2+||y||^2=1$, then $||F_n(x,y)||=||(2xy,||y||^2
-||x||^2)||=1$. 

Now by the Cayley-Dickson doubling process ([D]) define 
$$\a_{n+1}=\a_n\times\a_n$$ 
with
$$(a,b)(x,y)=(ax-\overline{y}b, ya+b\overline{x})
\,\,{\rm for}\,\,a,b,x\,\,{\rm  and}\,\, y\,\,{\rm  in}\,\, \a_n$$ 
and
$$\overline{x}=(\overline{x}_1,-x_2)\,\hbox{\rm if} \,\,x=(x_1,x_2)\,\,\hbox{\rm in}\,\,\a_{n-1}\times\a_{n-1}.$$
Thus, if $\a_0=\erre$ with $\overline{x}=x$ for $x$ real number,  then $\a_1=\ce,\a_2=\ache$ and $\a_3=\o$, which are normed algebras i.e.; $||xy||=||x||||y||$ for all $x,y$ in $\a_n$.

For $n\geq 4,\quad\a_n$ is no longer normed and also $\a_n$ has zero divisors
(see [K-Y] and [Mo$_1$]).

Let us define $X^\infty_n=\{(x,y)\in\a_n\times\a_n|F_n(x,y)=(0,0)\}$ and for $r$ nonnegative real number
$(x,y)\in X^r_n$  
if and only if $xy=0$ and $||x||=||y||=r$. It is clear that for $r>0$ and $s>0$ real numbers. $X^r_n$ is homeomorphic to $X^s_n$.
Let us define $X_n:=X^1_n$.

The set $X_n$ show up in some important problems in algebraic topology:

\noindent
(1) Cohen's approach to the Arf invariant one problem. (See $[C_1]$ and $[C_2]$).

\noindent
(2) Adem-Lam construction of normed and non-singular bilinear maps.(See [A] and [L]).

In this paper, we will show that for $n \geq4,\,\, X_n$ is related to some Stiefel manifolds; using the algebra structure in $\a_{n+1}$ we will construct the chain of inclusions
$$X_n\subset W_{2^{n-1}-1,2}\subset V_{2^n-2,2}\subset V_{2^n-1,2}$$
(see $\S$ 2 below) where $V_{m,2}$ and $W_{m,2}$ denote the real and complex Stiefel manifolds of 2-frames in $\erre^m$ and $\ce^m$ respectively.

In $\S$ 3 we show that we can attach  to every element in $W_{2^{n-1},2}$
in a canonical way,
 an eight 
dimensional vector subspace of $\a_{n+1}$ and that, only for the elements in $X_n$, such
vector subspace is isomorphic, as algebra, to $\a_3=\o$ (the octonions).

In $\S$ 4 we describe $X_n$ as a certain type of algebra monomorphisms from $\a_3=\o$ to $\a_{n+1}$ for $ n\geq 4$. 

This paper is a sequel of [Mo1] and we use freely the results of [Sch],we acknowledge with gratitude the hard work made by the reviwer.
%-------------------------
\vglue1cm
\noindent
{\bf \S 1. Pure and doubly pure elements in $\a_{n+1}.$} 
\vglue.5cm
\noindent

    Throughout this paper we use the following notational conventions:
    
 (1) Elements in $\a_n$ will be denoted by Latin characters  $a,b,c,\ldots, x,y,z.$
and elements in $\a_{n+1}$ will be denoted by Greek characters  $\alpha,\beta,\gamma,\ldots$.

 For example, 
$$\alpha=(a,b)\in\a_n\times\a_n.$$

 (2) When we need to represent  elements  in $\a_n$ as  elements  in $\a_{n-1}\times\a_{n-1}$ we use subscripts, for instance, $a=(a_1,a_2),\quad b=(b_1,b_2),$   and so on, with  $a_1,a_2,b_1,b_2$ in $\a_{n-1}.$
\vglue.2cm
\noindent

Now  $\{e_0,e_1,\ldots,e_{2^n-1}\}$ denotes the canonical basis in $\a_n.$ 
Then by the doubling process
$$\{(e_0,0),(e_1,0),\ldots,(e_{2^n-1},0),(0,e_0),\ldots,(0,e_{2^n-1})\}$$
is the canonical basis in $\a_{n+1}=\a_n\times\a_n$. By standard abuse of notation, we also denote $e_0=(e_0,0), e_1=(e_1,0),\ldots,e_{2^n-1}=(e_{2^n-1},0), e_{2^n}=(0,e_0),\ldots,e_{2^{n+1}-1}=(0,e_{2^n-1})$ in $\a_{n+1}.$ 

\vskip.5cm
\noindent

 For $\alpha=(a,b)\in\a_n\times\a_n=\a_{n+1}$ we denote $\widetilde{\alpha}=(-b,a)$ (the complexification of $\alpha$) so $\widetilde{e}_0=(0,e_0)$ and $\alpha\widetilde{e}_0=(a,b)(0,e_0)=(-b,a)=\widetilde{\alpha}$. Notice that $\widetilde{\widetilde{\alpha}}=-\alpha.$

\vskip.5cm
\noindent

\mbox{The {\bf trace}} on $\a_{n+1}$ is the linear map $t_{n+1}:\a_{n+1}\rightarrow\erre$ given by $t_{n+1}(\alpha)=\alpha+\overline{\alpha}=2$(real part of $\alpha$) 
so $t_{n+1}(\alpha)=t_n(a)$ when $\alpha=(a,b)\in\a_n\times\a_n.$
\vglue.5cm
\noindent
{\bf Definition:} $\alpha=(a,b)$ in $\a_{n+1}$ {\it is pure} if
$$t_{n+1}(\alpha)=t_n(a)=0.$$

$\alpha=(a,b)$ in $\a_{n+1}$ {\it is doubly pure} if it is pure and also $t_n(b)=0$;  i.e., $\widetilde{\alpha}$ is pure  in $\a_{n+1}.$

Also $2\langle a,b\rangle =t_n (a\overline{b})$ when $\langle -,-   \rangle $ is the inner product in $\erre^{2^n}$ (see [A]).

Note that for $a$ and $b$ pure elements $a\perp b$ if and only if $ab=-ba.$
\vglue.5cm
\noindent
{\bf Notation:} $\,_o\a_n=\{e_o\}^\perp\subset\a_n$ is the vector subspace consisting of pure elements in $\a_n$;  i.e., $\,_o\a_n={\rm Ker}(t_n)=\erre^{2^n-1}.$

$\widetilde{\a}_{n+1}=\,_o\a_n\times\,_o\a_n=\{e_0,\widetilde{e}_0\}^\perp=\erre^{2^{n+1}-2}$ is the vector subspace consisting of doubly pure elements in $\a_{n+1}.$
\vglue.5cm
\noindent
{\bf Lemma 1.1.} For $a$ and $b$ in $\widetilde{\a_n}$ we have that
\begin{itemize}
\item[1)] $a\widetilde{e}_0=\widetilde{a}$ and $\widetilde{e}_0a=-\widetilde{a}.$
\item[2)] $a\widetilde{a}=-||a||^2\widetilde{e}_0$ and $\widetilde{a} a=||a||^2\widetilde{e}_0$ so $a\perp \widetilde{a}.$
\item[3)] $\widetilde{a}b=-\widetilde{ab}$ with $a$ a pure element.
\item[4)] $a\perp b$ if and only if $\widetilde{a}b+\widetilde{b}a=0.$
\item[5)] $\widetilde{a}\perp b$ if and only if $ab=\widetilde{b}\widetilde{a}.$
\item[6)] $a\perp b$ and $\widetilde{a}\perp b$ if and only if $\widetilde{a}b=a\widetilde{b}.$
\end{itemize}
{\bf Proof:} Note that $a$ is pure if  and only if $\overline{a}=-a$ and if $a=(a_1, a_2)$ is doubly pure, then $\overline{a}_1=-a_1$ and $\overline{a}_2=-a_2.$
\begin{itemize}
\item[1)] $\widetilde{e}_0a=(0, e_0)(a_1, a_2)=(-\overline{a}_2, \overline{a}_1)=(a_2, -a_1)=-(-a_2,a_1)=-\widetilde{a}.$
\item[2)] $a\widetilde{a}=(a_1, a_2) (-a_2, a_1)=(-a_1a_2+a_1a_2, a_1^2+a^2_2)=(0,-||a||^2e_0)=-||a||^2\widetilde{e}_0.$

Similarly $\widetilde{a}a=(-a_2, a_1)(a_1, a_2)=(-a_2a_1+a_2a_1, -a^2_2-a^2_1)=||a||^2\widetilde{e}_0.$ 

Now, since $-2\langle \widetilde{a}, a\rangle =a\widetilde{a}+\widetilde{a}a=0$ we have $a\perp \widetilde{a}.$
\item[3)] $\widetilde{a}b=(-a_2, a_1)(b_1, b_2)=(-a_2b_1+b_2a_1, -b_2a_2 -a_1b_1)$.

{\it So $\widetilde{\widetilde{a}b}= (a_1b_1+b_2a_2, b_2a_1-a_2b_1)=(a_1, a_2)(b_1, b_2)=ab$ and  then 
$-\widetilde{a}b=\widetilde{ab}.$}

Notice that in this proof we only use that $\overline{a}_1=-a_1$; i.e., $a$ is pure and $b$ doubly pure.
\item[4)] $a\perp b\Leftrightarrow ab+ba=0\Leftrightarrow ab=-ba\Leftrightarrow\widetilde{ab}=-\widetilde{ba}.$

$\Leftrightarrow -\widetilde{a}b=\widetilde{b}a\Leftrightarrow \widetilde{a}b+\widetilde{b}a=0$ by (3).
\item[5)] $\widetilde{a}\perp b \Leftrightarrow \widetilde{\widetilde{a}}b+\widetilde{b}\widetilde{a}=0$ (by (4)) $\Leftrightarrow -ab+\widetilde{b}\widetilde{a}=0.$
\item[6)] If $\widetilde{a}\perp b$ and $a\perp b$, then by (3) and (4)
$
\widetilde{a}b=-\widetilde{ab}=\widetilde{ba}=-\widetilde{b}a=a\widetilde{b}.
$

Conversely, put $a=(a_1, a_2)$ and $b=(b_1, b_2)$ in $\a_{n-1}\times \a_{n-1}$ and define 

$c:=(a_1b_1+b_2a_2)$ and $d:=(b_2a_1-a_2b_1)$ in $\a_{n-1}.$

Then $a\widetilde{b}=(a_1, a_2)(-b_2, b_1)=(-a_1b_2+b_1a_2, b_1a_1+a_2b_2)$ so $a\widetilde{b}=(-\overline{d}, \overline{c}).$
\end{itemize}
Now $ab=(a_1, a_2)(b_1, b_2)=(a_1b_1+b_2a_2, b_2a_1-a_2b_1)=(c,d)$, so $\widetilde{ab}=(-d,c)$ and then $\widetilde{a}b=(d, -c)$. 

Thus, if $a\widetilde{b}=\widetilde{a}b$ then $\overline{c}=-c$ and $d=-\overline{d}.$
Then 
\begin{eqnarray*} 
t_n(ab)&=&t_{n-1}(c)=c+\overline{c}=0\;\;\hbox{\rm and}\;\;a\perp b\\
t_n(\widetilde{a}b)&=&t_{n-1}(d)=d+\overline{d}=0\;\;\hbox{\rm and}\;\;\widetilde{a}\perp b.
\end{eqnarray*}
\hfill Q.E.D

\vskip .5cm
\noindent
{\bf Corollary 1.2} For each $a\neq 0$ in $\widetilde{\a}_n$ the four dimensional vector subspace generated by $\{e_0, \widetilde{a}, a, \widetilde{e}_0\}$ is a copy of $\a_2 =\ache$ we denote it by $\ache_a.$
\vglue .25cm
\noindent
{\bf Proof:} We suppose that $||a||=1$, otherwise we take ${a\over ||a||}$. Construct the following multiplication table.

\vskip .5cm
\begin{center}
\begin{tabular}{c|cccc}
 & $e_0$ & $\widetilde{a}$ & $a$ & $\widetilde{e}_0$\\ \hline
$e_0$ & $e_0$ & $\widetilde{a}$ & $a$ & $\widetilde{e}_0$ \\
$\widetilde{a}$ & $\widetilde{a}$ & $-e_0$ & $+\widetilde{e}_0$ & $-a$\\
$a$ & $a$ & $-\widetilde{e}_0$ & $-e_0$ & $\widetilde{a}$\\
$\widetilde{e}_0$ & $\widetilde{e}_0$ & $a$ & $-\widetilde{a}$ & $-e_0$
\end{tabular}
\end{center}

By lemma 2.1. $a\widetilde{e}_0=\widetilde{a}$; $\widetilde{e}_0a=-\widetilde{a}$; $\widetilde{a}\widetilde{e}_0=\widetilde{\widetilde{a}}=-a;$ $\widetilde{e}_0 \widetilde{a}=-\widetilde{\widetilde{a}}=a$; $a\widetilde{a}=-\widetilde{e}_0$ and $\widetilde{a}a=\widetilde{e}_0.$

Identifying $e_0\leftrightarrow e_0$, $\widetilde{a}\leftrightarrow e_1,$ $a\leftrightarrow e_2$ and $\widetilde{e}_0\leftrightarrow e_3$ we have the multiplication table for $\a_2=\ache$

\hfill Q.E.D.
\vskip1cm
\noindent
{\bf \S 2.- The Stiefel manifold $V_{2^n-1,2}$ in $\a_{n+1}$ and a $T^2$-action.}

\vskip.5cm
Let $\langle a,b\rangle_n$ denote the standard inner product of $a$ and $b$ in $\a_n=\erre^{2^n}$. Now by [A] and [Mo$_1$]
$$2\langle a,b\rangle_n=(a\overline{b}+b\overline{a})=t_n(a\overline{b}).$$

It  is also  well known that, for $\alpha=(a,b)$ and $\chi=(x,y)$ in $\a_n\times\a_n=\a_{n+1}$ we have
$$\langle\alpha,\chi\rangle_{n+1}=\langle a,x\rangle_n+\langle\overline{b},\overline{y}\rangle_n.$$

In particular, if $\alpha$ and $\chi$ are doubly pure elements in $\a_{n+1}$ 
then  $y$ and $b$ are pure elements in $\a_n$, therefore
 
$$\langle\alpha,\chi\rangle_{n+1}=\langle a,x\rangle_n+\langle b,y\rangle_n.$$

\vskip.5cm
\noindent
{\bf Lemma 2.1} For $\alpha=(a,b)$ in $\a_{n+1}$ define  $\hat{\alpha}
:=(b,a)$.

For $\alpha\in\widetilde{\a}_{n+1}$ we have that

i) $\langle\alpha,\hat{\alpha}\rangle_{n+1}=0$ i.e. $\alpha\perp\hat{\alpha}$ in $\a_{n+1}$ if and only if $\langle a,b\rangle_n=0$, i.e. $a\perp b$ in $\a_n$.

ii) $\langle\widetilde{\alpha},\hat{\alpha}\rangle_{n+1}=0$, i.e. $\widetilde{\alpha}\perp\hat{\alpha}$ in $\a_{n+1}$ if and only if $||a||=||b||$ in $\a_n$.

\noindent
{\bf Proof:}  

i) $\langle\alpha,\hat{\alpha}\rangle_{n+1}=\langle(a,b),(b,a)\rangle_{n+1}=2\langle a,b\rangle_n.$

ii) $\langle\widetilde{\alpha},\hat{\alpha}\rangle_{n+1}=\langle(-b,a),(b,a)\rangle_{n+1}=-\langle b,b\rangle_n+\langle a,a\rangle_n=-||b||^2+||a||^2.$

\hfill Q.E.D.
\noindent

By \S 1, we know that for each $\alpha\neq 0$ in $\widetilde{\a}_{n+1}.
\,\,\ache_\alpha =$  Span $\{e_0,\widetilde{\alpha}, \alpha,\widetilde{e}_0\}$ is a copy of $\a_2$  and that if $\ache^\perp_\alpha$ denotes the orthogonal complement of $\ache_\alpha$, then $\a_{n+1}=\ache_\alpha\oplus\ache^\perp_\alpha$.

 Since $\alpha$ is doubly pure, $\hat{\alpha}$ is also doubly pure; i.e., $\hat{\alpha}\in\{e_0,\widetilde{e}_0\}^\perp$ and if $\hat{\alpha}\perp\alpha$ and $\hat{\alpha}\perp\widetilde{\alpha},$ then $\hat{\alpha}\in\ache^\perp_\alpha$. Now $S^{\sqrt{2}}(\widetilde{\a}_{n+1})=S^{2^{n+1}-3}$ denotes the sphere of radius $\sqrt{2}$ inside of $\widetilde{\a}_{n+1}.$

 Thus, we have a description of the real Stiefel manifold of 2- orthonormal frames in $\erre^{2^n-1}$ as follows:
$$V_{2^n-1,2}=\{(a,b)\in\,_o\a_n\times\,_o\a_n=\widetilde{\a}_{n+1}:||a||=||b||=1, a\perp b\}$$
and
$$V_{2^n-1,2}=\{\alpha\in S^{\sqrt{2}}(\widetilde{\a}_{n+1}) : \hat{\alpha}\in\ache^\perp_\alpha\}.$$

\noindent
{\bf Lemma 2.2.} If $r$ and $s$ are in $\erre$ with $r^2+s^2=1$ and $(a,b)\in V_{2^n-1,2}$ then $(ra-sb,sa+rb)\in V_{2^n-1,2}.$

\noindent
{\bf Proof:} Suppose that 
$||a||=||b||=1$ and $a\perp b$ in $\a_n$. Then $||ra-sb||^2=r^2||a||^2+s^2||b||^2-2rs\langle a,b\rangle_n$ and $||sa+rb||^2=s^2||a||^2+r^2||b||^2+2rs\langle a,b\rangle_n$, so $||ra-sb||^2=||sa+rb||^2=r^2+s^2=1.$
\begin{eqnarray*}
\langle ra-sb,sa+rb\rangle_n&=&rs\langle a,a\rangle_n-sr\langle b,b\rangle_n-s^2\langle b,a\rangle_n+r^2\langle a,b\rangle_n\\
&=&rs||a||^2-rs||b||^2+0=rs-sr=0.
\end{eqnarray*}
\hfill Q.E.D.
\vglue.5cm
\noindent
{\bf Corollary 2.3.} $S^1\times V_{2^n-1,2}\stackrel{{\rm a}}{\rightarrow}V_{2^n-1,2}$ given by 
$$((r,s),\alpha)\mapsto r\alpha+s\widetilde{\alpha}=(ra-sb,sa+rb)$$ 
defines a smooth, free $S^1$-action on $V_{2^n-1,2}.$
\vskip.5cm
\noindent
{\bf Proof:} Clearly $(1,0)\cdot\alpha=\alpha$ and 

$(r,s)[(q,t)\cdot\alpha]=((r,s)(q,t))\cdot\alpha=(rq-st,rt+sq)\cdot\alpha$ so ${\rm a}$ defines an action. It is a smooth action because it is a restriction of a linear action of $GL_2(\erre)$ on $\widetilde{\a}_{n+1}=\erre^{2^{n+1}-2}.$

Finally, ${\rm a}$ is a free action: if $r\alpha+s\widetilde{\alpha}=\alpha$ then $r=1$ and $s=0$, because $\alpha\perp\widetilde{\alpha}$ in $\widetilde{\a}_{n+1}.$

\hfill Q.E.D.

Now, we identify $V_{2^n-2,2}$,  the real Stiefel manifold of 2-orthonormal frames on $\erre^{2^n-2}$, as a submanifold of $V_{2^n-1,2}$ as follows:
$$V_{2^n-2,2}=\{(a,b)\in V_{2^n-1,2}|(a,b)\in\widetilde{\a}_n\times\widetilde{\a}_n\};$$
i.e., $(a,b)\in V_{2^n-1,2}$ belongs to $V_{2^n-2,2}$ whenever $a$ and $b$ are doubly pure elements in $\a_n$ and we have the known fibration [Wh]
$$
\begin{array}{lllll}
S^{2^n-4}&\rightarrow&V_{2^n-2,2}&\rightarrow&S(\widetilde{\a}_n)=S^{2^n-3}\\
&&(a,b)&\mapsto&b
\end{array}
$$
thus, $V_{2^n-2,2}$ has dimension $2^n-3+2^n-4=2^{n+1}-7.$
 
Since $(ra-sb)$ and $(sa+rb)$ are doubly pure elements in $\a_n$ when $a$ and $b$ are doubly pure elements, we have that 

{\bf $V_{2^n-2,2}$ is a $S^1$-invariant submanifold of $V_{2^n-1,2};$} 

i.e., if $\alpha\in V_{2^n-2,2}$ then $(r,s)\cdot\alpha\in V_{2^n-2,2}$ for all $(r,s)\in S^1.$

\vskip.3cm
\noindent
We note that $\a_{n+1}$ becomes a complex vector space by defining $i\dot \alpha$=$\widetilde{\alpha}$ thus as a complex vector space
$$\widetilde{\a}_{n+1}= \,_0\!\a_n\times\,_0\!\a_n \cong\ce \otimes_{\errita}  \,_0\!\a_n 
$$
The isomorphism takes $1\otimes x$ to $(x,0)$ and $i\otimes y$ to $(0,y)$ and $S^1$ (the set of modulo 1 complex numbers) acts naturally by multiplication on $\ce$, hence on $\widetilde{\a}_{n+1}.$

Now, we identify the complex Stiefel manifold $W_{2^{n-1}-1,2}$  of 2-orthonormal frames in $\ce^{2^{n-1}-1}$ as a submanifold of $V_{2^n-2,2}$ in terms of the Cayley-Dickson algebra $\a_{n+1}$ for $n\geq 3.$

It is known   that for $a,b$ and $x$ in $\a_n, \quad\langle ax,b\rangle_n=\langle a,b\overline{x}\rangle_n$.(See [A]). Thus if $x$ is a pure element, i.e., $\overline{x}=-x$ then $\langle ax,b\rangle_n=-\langle a,bx\rangle_n.$ That is, right multiplication by a pure non-zero element is a skew-symmetric linear map. In particular $\langle\widetilde{a},b\rangle_n=-\langle a,\widetilde{b}\rangle_n.$

\vskip.3cm
\noindent
{\bf Proposition 2.4.} 
For $n\geq 3,$ the map ${\cal H}_n:\widetilde{\a}_n\times\widetilde{\a}_n\rightarrow\ce$  given by 

$${\cal H}_n(a,b)=2\langle a,b\rangle_n-2i\langle\widetilde{a},b\rangle_n$$ 
defines a Hermitian inner product in $\widetilde{\a}_n.$
\vskip.3cm
\noindent
{\bf Proof:} Clearly ${\cal H}_n$ is $\erre$-linear and
\begin{eqnarray*}
\overline{{\cal H}_n(a,b)}&=&2\langle a,b\rangle_n+2i\langle\widetilde{a},b\rangle_n\\
&=&2\langle a,b\rangle_n-2i\langle a,\widetilde{b}\rangle_n\\
&=&{\cal H}_n(b,a).
\end{eqnarray*}

On the other hand,
\begin{eqnarray*}
{\cal H}_n(\widetilde{a},b)&=&2\langle\widetilde{a},b\rangle_n-2i
\langle\widetilde{\widetilde{a}},b\rangle_n\\
&=&2\langle\widetilde{a},b\rangle_n+2i\langle a,b\rangle_n\\
&=&2i\langle a,b\rangle_n-2i^2\langle\widetilde{a},b\rangle_n\\
&=&i{\cal H}_n(a,b).
\end{eqnarray*} 

\hfill Q.E.D.

\vskip.3cm

\noindent
{\bf Proposition 2.5.} For $n\geq 3$
$$W_{2^{n-1}-1,2}=\{(a,b)\in V_{2^n-2,2} |  b\in\ache^\perp_a\}.$$
{\bf Proof:} First of all, we observe that $b\in\ache^\perp_a$ for $a$ and $b$ in $\widetilde{\a}_n$ if and only if $b\perp a$ and $b\perp\widetilde{a}$, i.e.
 ${\cal H}_n(a,b)=0$. If $||a||=||b||=1$ and ${\cal H}_n(a,b)=0$ then $(a,b)\in W_{m,2}$, where $m=\frac{1}{2}(2^n-2)=2^{n-1}-1.$

\hfill Q.E.D.

\vskip.3cm

\noindent
{\bf Proposition 2.6.} $W_{2^{n-1}-1,2}$ is $S^1$-invariant.
\vskip.3cm       
\noindent
{\bf Proof:} Suppose $(a,b)\in\widetilde{\a}_n\times\widetilde{\a}_n$ with $||a||=||b||=1$ and $b\in\ache^\perp_a.$ From this, we have $b\perp a$, and $\widetilde{b}\perp a$ (equivalently  $\widetilde{a}\perp b).$

Now $r(a,b)+s(-b,a)=(ra-sb,rb+sa)$ and we know that  $(ra-sb)\perp (rb+sa).$ 

To finish, we need to show that $\widetilde{(ra-sb)}\perp(rb+sa).$
\begin{eqnarray*}
\langle\widetilde{ra-sb},rb+sa\rangle_n&=&\langle r\widetilde{a}-s\widetilde{b},rb+sa\rangle_n\\
&=&r^2\langle\widetilde{a},b\rangle_n-s^2\langle\widetilde{b},a\rangle_n+rs\langle\widetilde{a},a\rangle_n-rs\langle\widetilde{b},a\rangle_n\\
&=&0,
\end{eqnarray*}
therefore $(ra-sb)\in\ache_{rb+sa}.$

\hfill Q.E.D.
\vskip.3cm
\noindent

Note that we have a fibration
$$
\begin{array}{lllll}
S^{2^n-5}&\rightarrow&W_{2^{n-1}-1,2}&\stackrel{\pi}{\rightarrow}
&S(\widetilde{\a}_n)=S^{2^n-3}\\
&&(a,b)&\mapsto&b
\end{array}
$$

$\pi^{-1}(b)=S(\ache^\perp_b)=S^{2^n-5}$ since dim $\ache^\perp_b=2^n-4.$
\vglue.2cm

Thus dim $W_{2^{n-1}-1,2}=2^n-5+2^n-3=2^{n+1}-8.$
\vglue.5cm

In [Mo$_1$] it is shown  that for $a$ and $b$ in $\a_n$ with $n\geq 4$ and $||a||=||b||=1,$ then

if $ab=0$ then

i) $(a,b)\in\widetilde{\a}_n\times\widetilde{\a}_n.$

ii) $b\in \ache^\perp_a$ (or equivalently $a\in \ache^\perp_b$). 

Thus
$$X_n:=\{(a,b)\in \a_n\times\a_n: ||a||=||b||=1\quad{\rm and}\quad ab=0\}$$
is a subset of $W_{2^{n-1}-1,2}.$

Thus we have a chain of inclusions for $n\geq 3,$

$$X_n\subset W_{2^{n-1}-1,2}\subset V_{2^n-2,2}\subset V_{2^n-1,2}.$$
\vglue.5cm
Now we show that 
$X_n$ and $W_{2^{n-1}-1,2}$ admit a $T:=S^1\times S^1$ action.
\vglue.5cm
\noindent
{\bf Lemma 2.7} For $(a,b)\in V_{2^n-2,2}$ and $r,s,q,p$ in $\erre$ with $r^2+s^2=1$ and $p^2+q^2=1$ define
$$(a,b)\stackrel{\tau}{\mapsto}(ra+s\tilde{a},pb+q\tilde{b}).$$

i) If $(a,b)\in W_{2^{n-1}-1,2}$ then $(ra+s\tilde{a},pb+q\tilde{b})\in W_{2^{n-1}-1,2}$.

ii) If $(a,b)\in X_n$ then $(ra+s\tilde{a},pb+q\tilde{b})\in X_n$.

iii) $\tau $ defines a free $T$-action on $W_{2^{n-1}-1,2}$ and $X_n$ respectively.
\vglue.5cm
\noindent
{\bf Proof:} By direct calculations. If $a\perp b$ and $\tilde{a}\perp b$ then 
\begin{eqnarray*}
\langle ra+s\tilde{a},pb+q\tilde{b}\rangle_n&=&rp\langle a,b\rangle_n+sq\langle\tilde{a},\tilde{b}\rangle_n+rq\langle a,\tilde{b}\rangle_n+sp\langle\tilde{a},b\rangle_n\\
&=&0+0+0+0\\
&=&0
\end{eqnarray*}
similarly.
\begin{eqnarray*}
\langle ra +s\tilde{b},(pb+q\tilde{b})\tilde{e_0}\rangle_n&=&\langle ra+s\tilde{a},p\tilde{b}-qb\rangle_n\\
&=&rp\langle a,\tilde{b}\rangle_n+sp\langle\tilde{a},\tilde{b}\rangle_n-rq\langle a,b\rangle_n-sq\langle\tilde{a},b\rangle_n\\
&=&0+0+0+0\\
&=&0.
\end{eqnarray*}
If $ab=0$ then
\begin{eqnarray*}
(ra+s\tilde{a})(pb+q\tilde{b})&=&
rp(ab)+sq\tilde{a}\tilde{b}+sp\tilde{a}b+rqa\tilde{b}\\
&=&0.
\end{eqnarray*}
Also
\begin{eqnarray*}
||ra+s\tilde{a}||^2&=&r^2||a||^2+s^2||\tilde{a}||^2=(r^2+s^2)||a||^2=1\quad{\rm and}\\
||pb+q\tilde{b}||^2&=&p^2||b||^2+q^2||\tilde{b}||^2=(p^2+q^2)||b||^2=1
\end{eqnarray*}
Therefore we proved (i) and (ii).

Finally $(ra+s\tilde{a},pb+q\tilde{b})=(a,b)$ if and only if $r=1, s=0,  p=1$ and $q=0$.

Clearly this action is  smooth and  free  (see Corollary 2.5).

\hfill Q.E.D.

\vglue1cm
\noindent
{\bf $\S$ 3. $X_n$, Octonions and a $S^3$ action.}
\vglue.5cm
In this section, we show that we can attach to every element in $X_n$ a copy of $\a_3=\o$, the octonions inside of $\a_{n+1}$ for $n>3$.This allows us to identify $X_n$ with a subset of algebra monomorphisms of $\a_3$ into $\a_{n+1},$ which is our main goal in $\S$ 4.

We recall some notation from $\S$ 1.

Let $e_0\in\a_{n-1}$ be the unit  so $(e_0,0)=e_0$ is the unit in $\a_n$ and $\widetilde{e}_0=(0,e_0)$ in $\a_n$.

For $\widetilde{e}_0$ in $\a_n$ we denote 
$\varepsilon=(\widetilde{e}_0,0)$ in $\a_{n+1}$.

For example for $n=4\,,\,\widetilde{e}_0=e_8$ in $\a_4$, and $\varepsilon=(e_8,0)$ in $\a_5$

In general $\widetilde{e}_0=e_{2^n}$ in $\a_{n+1}$ and $\varepsilon=e_{2^{n-1}}$ in $\a_{n+1}$. Since $\varepsilon$ is a doubly pure element of norm one, we have that $\ache_\varepsilon={\rm Span}\{e_0,\widetilde{\varepsilon},\varepsilon,\widetilde{e}_0\}\subset\a_{n+1}$ is a copy of $\a_2$ and a direct sum decomposition $\a_{n+1}=\ache_\varepsilon\oplus\ache^\perp_\varepsilon.$ 

By definition $\alpha=(a,b)\in\a_n\times\a_n=\a_{n+1}$ is doubly pure in $\a_{n+1}$ with doubly pure entries in $\a_n$ if and only if $\alpha \in\ache^\perp_\varepsilon$.

In section $\S$ 2 we constructed the chain
$$X_n\subset W_{2^{n-1}-1,2}\subset V_{2^n-2,2}\subset V_{2^n-1,2}\subset\a_{n+1}$$
for $n\geq 3$ with $X_3=\Phi$, the empty set.
\vglue.5cm
Therefore, by definition, $V_{2^n-2,2}=V_{2^n-1,2}\cap\ache^\perp_\varepsilon$.
\vglue.5cm
\noindent
{\bf Lemma 3.1.} For 
$\alpha\in\ache^\perp_\varepsilon\subset\a_{n+1}$ with 
$\alpha=(a,b)\in\widetilde{\a}_n \times\widetilde{\a}_n$,

1) $(\alpha\varepsilon)\in\ache^\perp_\varepsilon$ and 
$\alpha\varepsilon=(\widetilde{a},-\widetilde{b})$.

2) $\alpha\widetilde{\varepsilon}\in\ache^\perp_\varepsilon$ and 
$\alpha\widetilde{\varepsilon}=\widetilde{\alpha}\varepsilon=-\widetilde{\alpha\varepsilon}=(-\tilde{b},-\tilde{a}).$
\vglue.5cm
\noindent
{\bf Proof:} By direct calculation
$$\alpha\varepsilon=(a,b)(\widetilde{e}_0,0)=(a\tilde{e}_0,-b\tilde{e}_0)=(\tilde{a},-\tilde{b})\in\tilde{\a}_n\times\tilde{\a}_n=\ache^\perp_\varepsilon$$
and 
$$\alpha\tilde{\varepsilon}=(a,b)(0,\widetilde{e}_0)=
(\tilde{e}_0b,\tilde{e}_0a)=(-\tilde{b},-\tilde{a})$$
by Lemma 1.1 (1).

Finally using  Lemma 1.1 (6) and (3) respectively $\alpha\tilde{\varepsilon}=\tilde{\alpha}\varepsilon=-\widetilde{\alpha\varepsilon}\in\ache^\perp_\varepsilon.$
\vglue.5cm

\hfill Q.E.D.
\vglue.5cm
\noindent
{\bf Corollary 3.2.} For a non-zero $\alpha$ in 
$\ache^\perp_\varepsilon\subset\a_{n+1}$ and $ n\geq 3$,
$$\o_\alpha := Span \{e_0,\tilde{\varepsilon},\varepsilon,\widetilde{e}_0,\tilde{\alpha},\alpha\varepsilon,\tilde{\varepsilon}\alpha,\alpha\}\subset\a_{n+1}$$
is an 8-dimensional vector subspace of $\a_{n+1}=\erre^{2^{n+1}}$.
\vglue.5cm
\noindent
{\bf Proof:} By definition $\{e_0,\tilde{\varepsilon},\varepsilon,\tilde{e}_0\},\{e_0,\tilde{\alpha},\alpha,\tilde{e}_0\},\{\varepsilon,\alpha,\tilde{\varepsilon},\tilde{\alpha}\}$ are an orthogonal set of vectors and $\alpha\varepsilon\in\ache^\perp_\varepsilon\cap\ache^\perp_\alpha$. Also by Lemma 3.1. $\tilde{\varepsilon}\alpha=-\alpha\tilde{\varepsilon}\in\ache^\perp_\varepsilon\cap\ache^\perp_\alpha$. Thus $\{e_0,\tilde{\varepsilon},\varepsilon,\tilde{e}_0,\tilde{\alpha},\alpha\varepsilon,\tilde{\varepsilon}\alpha,\alpha\}$ is an orthogonal set of vectors in $\a_{n+1}$.

\hfill Q.E.D.
\vglue.5cm
\noindent
{\bf Remark:} In particular for $\alpha\in V_{2^n-2,2}$ we have that $\o_\alpha\cong\erre^8\subset\a_{n+1}$ and $\o_\alpha\oplus\o^\perp_\alpha=\a_{n+1}$.
\vglue.5cm
\noindent
{\bf Lemma 3.3.} For $\alpha\in V_{2^n-2,2}$.

$\alpha\in W_{2^{n-1}-1,2}$ if and only if $\hat{\alpha}\in\o^\perp_\alpha$.
\vglue.5cm
\noindent
{\bf Proof:} Recall that by definition $\hat{\alpha}=(b,a)$ if $\alpha=(a,b)$ so 

$\hat{\alpha}\in({\rm Span}((\{e_0,\tilde{\varepsilon},\varepsilon,\tilde{e}_0,\alpha,\tilde{\alpha}\}))^\perp$ (see Lemma 2.1. above).

Now 
\begin{eqnarray*}
\langle\hat{\alpha},\tilde{\varepsilon}\alpha\rangle_{n+1}&=&\langle(b,a),(\tilde{b},\tilde{a})\rangle_{n+1}=\langle b,\tilde{b}\rangle_n+\langle a,\tilde{a}\rangle_n=0\\
\langle\hat{\alpha},\alpha\varepsilon\rangle_{n+1}&=&
\langle(b,a),(\tilde{a},-\tilde{b})\rangle_{n+1}=
\langle b,\tilde{a}\rangle_n-\langle a,\tilde{b}\rangle_n=2\langle b,\tilde{a}\rangle_n,
\end{eqnarray*}
so $\hat{\alpha}\perp(\alpha\varepsilon)$ in $\a_{n+1}$ if and only if $\tilde{a}\perp b$ in $\a_n$ i.e. $b\in\ache^\perp_a$.

Therefore $\hat{\alpha}\in\o^\perp_\alpha$ if and only if $b\in\ache^\perp_a$.

\hfill Q.E.D.
\vglue.5cm
Thus $W_{2^{n-1}-1,2}=\{\alpha\in V_{2^n-2,2}|\hat{\alpha}\in\o^\perp_\alpha\}$.\vglue.5cm
\noindent
{\bf Theorem 3.4} For $\alpha\in W_{2^{n-1}-1,2}$ and $n\geq 4$,  the following statements are equivalent.

i) $\alpha\in X_n$

ii) $\alpha$ alternate with $\varepsilon$ i.e., $(\alpha,\alpha,\varepsilon)=0$

iii) The vector subspace of $\a_{n+1}$
$$V(\alpha;\varepsilon):= {\rm Span} \{e_0,\alpha,\varepsilon,\alpha\varepsilon\}.$$
is multiplicatively closed and isomorphic to $\a_2=\ache$.

iv) $\o_\alpha$ is multiplicatively closed and isomorphic to $\a_3=\o$.

v) $\hat{\alpha}\in{\rm Ker} L_\alpha\subset\o^\perp_\alpha$, where $L_\alpha$ is left multiplication by $\alpha$.
\vglue.5cm
\noindent
{\bf Proof.} First of all, we calculate
\begin{eqnarray*}
\alpha(\alpha\varepsilon)&=&(a,b)[(a,b)(\tilde{e},0)]=(a,b)(\tilde{a},-\tilde{b})=(a\tilde{a}-\tilde{b}b,-\tilde{b}a-b\tilde{a})\\
&=&(-||a||^2\tilde{e}_0-||b||^2\tilde{e}_0,-\tilde{b}a-\tilde{b}a)
\,\,\hbox{(\rm by Lemma 1.1 (2) and (5))}\\ 
&=&-||\alpha||^2\varepsilon+2(0,\tilde{ba}).
\end{eqnarray*}
Therefore $\alpha(\alpha\varepsilon)=\alpha^2\varepsilon=-||\alpha||^2\varepsilon$ if and only if $ba=0$, i.e. $\alpha\in X_n$ and we have (i)$\Leftrightarrow$ (ii).

Clearly if $\alpha\in W_{2^{n-1}-1,2}$ then $\{e_0,\alpha,\varepsilon,\alpha\varepsilon\}$ form an orthonormal set of vectors so ${\rm \dim}_\erre(V(\alpha;\varepsilon))=4$ so $-||\alpha||^2=\alpha^2=(\alpha\varepsilon)^2$ and $\alpha(\alpha\varepsilon)=-||\alpha||^2\varepsilon$ if and only if $V(\alpha;\varepsilon)=\ache$, and we prove (ii)$\Leftrightarrow$(iii). 

To prove (iii)$\Leftrightarrow$(iv) we stablish the following correspondence between the canonical basis in $\a_3$ and the orthonormal basis of $\o_\alpha.$
$$e_1\rightarrow\tilde{\varepsilon}; e_2\mapsto\varepsilon;e_3\rightarrow\tilde{e}_0; ||\alpha||e_4\rightarrow\tilde{\alpha};||\alpha||e_5\rightarrow\alpha\varepsilon; ||\alpha||e_6\rightarrow\tilde{\varepsilon}\alpha; ||\alpha||e_7\rightarrow\alpha$$
Using ii) it is a routine calculation to see that this correspondence define an algebra isomorphism.
(See also Lemma 4.4 (1) below).

Finally by Lemma 3.3. we know $\hat{\alpha}\in\o^\perp_\alpha$ and
$$\alpha\hat{\alpha}=(a,b)(b,a)=(ab+ab,a^2-b^2)=(2ab,||b||^2-||a||^2)$$ 
is the Hopf construction. So 
$\alpha\hat{\alpha}=0$ in $\a_{n+1}$ if and only if  $\alpha\in X_n.$
Recall that $\a_3\cong\o_\alpha$ 
 does not admit zero divisors.

\hfill Q.E.D.
\vglue.5cm
\noindent
{\bf Theorem 3.5.} $\ache^\perp_\varepsilon$ admits a left $\ache_\varepsilon$-module structure for $n\geq 3$.
\vglue.5cm
\noindent
{\bf Proof:} For $\alpha=(a,b)$ in $\ache^\perp_\varepsilon=\widetilde{\a}_n\times\widetilde{\a}_n$ and $u=re_0+s\widetilde{\varepsilon}+q\varepsilon+p\widetilde{e}_0$ with $r,s,q$ and $p$ in $\erre$.

Define $$u\cdot\alpha=\alpha u=r\alpha+s\alpha
\widetilde{\varepsilon}+q\alpha\varepsilon+p\widetilde{\alpha}.$$

Trivially 
$\widetilde{\alpha}\in\ache^\perp_\varepsilon$ and $(\alpha
\widetilde{\varepsilon})$ and $(\alpha\varepsilon)$ are in $\ache_\varepsilon^\perp$ by Lemma 3.1.(2) and (1) respectively.

Since 
$\widetilde{\varepsilon},\varepsilon$ and 
$\widetilde{e}_0$ are alternative elements in 
$\a_{n+1}$ (actually they belong to the canonical basis) we have that 
$\widetilde{\varepsilon}\cdot\alpha=(\alpha\widetilde{\varepsilon})
\widetilde{\varepsilon}=\alpha(\widetilde{\varepsilon})^2=-||\alpha||^2e_0=
\widetilde{\varepsilon}\cdot\alpha$ and similarly 
$\varepsilon\cdot(\varepsilon\cdot\alpha)=\varepsilon^2\cdot\alpha$ and 
$\widetilde{e}_0\cdot(\widetilde{e}_0\cdot\alpha)=
\widetilde{e}^2_0\cdot\alpha$.

Now $\varepsilon\cdot(\widetilde{e}_0\cdot\alpha)=\varepsilon\cdot(\alpha\widetilde{e}_0)=\varepsilon\cdot\widetilde{\alpha}=\widetilde{\alpha}\varepsilon$ and
$(\varepsilon\widetilde{e}_0)\cdot\alpha=\widetilde{\varepsilon}\cdot\alpha=\alpha\widetilde{\varepsilon}=\widetilde{\alpha}\varepsilon$ by Lemma 3.1 (2).

Similarly 
\begin{eqnarray*}
\widetilde{\varepsilon}\cdot(\widetilde{e}_0\cdot\alpha)&=&(\widetilde{\varepsilon}\widetilde{e}_0)\cdot\alpha=\varepsilon\alpha\\
\widetilde{e}_0\cdot(\widetilde{\varepsilon}\cdot\alpha)&=&(\widetilde{e}_0\widetilde{\varepsilon})\cdot\alpha=\alpha\varepsilon\\
\varepsilon\cdot(\widetilde{e}_0\cdot\alpha)&=&(\varepsilon\widetilde{e}_0)\cdot\alpha=\widetilde{\alpha}\varepsilon\\
\widetilde{e}_0\cdot(\varepsilon\cdot\alpha)&=&(\widetilde{e}_0\varepsilon)\cdot\alpha=-\alpha\widetilde{\varepsilon}.
\end{eqnarray*}
Finally $\widetilde{\varepsilon}
\cdot(\varepsilon\cdot\alpha)=
\widetilde{\varepsilon}\cdot(\alpha\varepsilon)=(\alpha\varepsilon)
\widetilde{\varepsilon}=(
\widetilde{\alpha\varepsilon})\varepsilon=-
\widetilde{(\alpha\varepsilon)\varepsilon}=
\widetilde{-\alpha}=(
\widetilde{\varepsilon}\varepsilon)\cdot\alpha$
and 
$\varepsilon\cdot(
\widetilde{\varepsilon}\cdot\alpha)=
\varepsilon\cdot(\alpha
\widetilde{\varepsilon})=(\alpha
\widetilde{\varepsilon})\varepsilon=(
\widetilde{\alpha}\varepsilon)\varepsilon=
\widetilde{\alpha}=\varepsilon\cdot(\widetilde{\varepsilon}\alpha)$.

By Lemma 3.1 and Lemma 1.1. and we are done.

\hfill Q.E.D.

\vglue.5cm
Now we define a $S^3$ action on $X_n$.

Consider the unit sphere inside of $\ache_\varepsilon\subset\a_{n+1}$.

$$S^3=S(\ache_\varepsilon)=\{re_0+s\tilde{\varepsilon}+q\varepsilon+p
\tilde{e}_0|r^2+s^2+q^2+p^2=1\}.$$

For $\alpha\in\ache^\perp_\varepsilon$ with $\alpha=(a,b)\in\tilde{\a}_n\times\tilde{\a}_n$ define 
$\ache^\perp_\varepsilon\times S^3\rightarrow\ache^\perp_\varepsilon$ by
\begin{eqnarray*}
\alpha(re_0+s\tilde{\varepsilon}+q\varepsilon+p\tilde{e}_0)&=&r\alpha+s\alpha\tilde{\varepsilon}+q\alpha\varepsilon+p\alpha\tilde{e}_0=r\alpha+s\tilde{\alpha}\varepsilon+q\alpha\varepsilon+p\tilde{\alpha}\\
&=&r(a,b)+s(-\tilde{b},-\tilde{a})+q(\tilde{a},-\tilde{b})+p(-b,a)\\
&=&(ra-s\tilde{b}+q\tilde{a}-pb,rb-s\tilde{a}-q\tilde{b}+pa).
\end{eqnarray*}
By definition this is a group action which is smooth and free of fixed points.
\vglue.5cm
\noindent
{\bf Corollary 3.6} The above action of $S^3=S(\ache_\varepsilon)$ on $\ache^\perp_\varepsilon$ is a group action which is smooth, orthogonal and free of fixed points.
\vglue.5cm
\noindent
{\bf Proof:} By Theorem 3.5 this is a smooth group action because it is a restriction of a linear action.  Since right multiplication by $\widetilde{e}_0,\varepsilon$ and $\widetilde{\varepsilon}$ are orthogonal linear transformations, we have that the action is orthogonal.

Finally this action is free of fixed points because $\{e_0,\widetilde{\varepsilon},\varepsilon,\varepsilon\}$ is an orthonormal basis, so $\alpha(re_0+s\widetilde{\varepsilon}+q\varepsilon+p\widetilde{e}_0)=\alpha$ if and only if $r=1, s=q=p=0$.

\hfill Q.E.D.
\vglue.5cm
\noindent
{\bf Theorem 3.8.}

 i) The subsets $X_n$ and $W_{2^{n-1}-1,2}$ of $\ache^\perp_\varepsilon$ are $S^3$-equivariant.

ii) For $\alpha$ and $\beta$ in $W_{2^{n-1}-1,2,}$

$\o_\alpha=\o_\beta,$ as vector spaces, if and only if $\alpha$ and $\beta$ lie in the same $S^3$-orbit.
\vglue.5cm
\noindent
{\bf Proof.} For $\alpha\in\ache^\perp_\varepsilon$ with $\alpha=(a,b)\in\tilde{\a}_n\times\tilde{\a}_n$ and $r,s,q$ and $p$ in $\erre$ with $r^2+s^2+q^2+p^2=1$ we have that
$$\alpha(re_0+s\tilde{\varepsilon}+q\varepsilon+s\tilde{e}_0)=(ra-s\tilde{b}+q\tilde{a}-pb,rb-s\tilde{a}-q\tilde{b}+pa).$$
Suppose that $\alpha\in W_{2^{n-1}-1,2}$ then $\langle\tilde{a},b\rangle_n=-\langle a,\tilde{b}\rangle_n=0;\langle a,b\rangle_n=\langle\tilde{a},\tilde{b}\rangle_n=0$  and by definition $\langle a,\tilde{a}\rangle_n=\langle b,\tilde{b}\rangle_n=0$ with $||a||=||\tilde{a}||=||\tilde{b}||=||b||=1$ so $\langle ra-s\tilde{b}+q\tilde{a}-pb,rb-s\tilde{a}-q\tilde{b}+pa\rangle= p||a||^2+sq||b||^2-qs||\tilde{a}||^2-pr||b||^2=0$ and $(\alpha(re_0+s\tilde{\varepsilon}+q\varepsilon +s\tilde{e}_0))\in V_{2^n-2,2}$.

Similarly $\langle ra-s\tilde{b}+q\tilde{a}-pb,r\tilde{b}+sa+qb+p\tilde{a}\rangle=0$ and 
$$\alpha(re_0+s\tilde{\varepsilon}+q\varepsilon +s\tilde{e}_0)\in W_{2^{n-1}-1,2}.$$

A direct calculation shows that if $ab=0$ then
$$ab=\tilde{a}b=a\tilde{b}=\tilde{a}\tilde{b}=0$$

$$(ra-s\tilde{b}+q\tilde{a}-pb)(rb-s\tilde{a}-q\tilde{b}+pa)=$$

$$-rsa\tilde{a}+rpa^2-sr\tilde{b}b+sq\tilde{b}^2-qs\tilde{a}^2+qpaa-prb^2+pqb\tilde{b}=$$

$$-rs(-||a||^2\tilde{e}_0+||b||^2\tilde{e}_0)+pq(-||a||^2\tilde{e}_0+||b||^2\tilde{e}_0)+rp(a^2-b^2)+sq(\tilde{b}^2-\tilde{a}^2)=0$$
because $||a||^2=||b||^2=1$ and $a^2=b^2=-e_0$, so we have (i). 

To prove ii) we notice that $\o_\alpha=\ache_\varepsilon\oplus {\rm Span}\{\tilde{\alpha},\alpha\varepsilon,\tilde{\varepsilon}\alpha,\alpha\}.$Then if $\beta=r\alpha+s\tilde{\alpha}\varepsilon+q\alpha\varepsilon+p\tilde{\alpha}$ (recall that $\tilde{\alpha}\varepsilon=-\tilde{\varepsilon}\alpha$ by Lemma 4.1) and $||\beta||=1$ then $r^2+s^2+q^2+p^2=1$ and $\alpha\equiv\beta$ mod $S^3$ if and only if $\o_\beta\subset\o_\alpha$ but 
$$\dim \o_\beta= \dim \o_\alpha=8\,\,{\rm  and}\,\, \o_\beta=\o_\alpha .$$

\hfill Q.E.D.

{\bf Remark:} Notice that $T=S^1\times S^1$, as in Lemma 2.7. and 
$S^3=S(\ache_\varepsilon)$ intersect on a copy of $S^1$.

Suppose that 
$r^2+s^2+p^2+q^2=1$ and $u^2+v^2=1\quad t^2+m^2=1$ in $\erre$.

If $(ra-s\tilde{b}+q\tilde{a}-pb,rb-s\tilde{a}-q\tilde{b}+pa)=(ua+v\tilde{a},tb+m\tilde{b})$
then $r=u$,

$q=v,s=0,p=0,r=t,-q=m$ so

$$S(\ache_\varepsilon)\cap T=S^1=\{(r,-q)|r^2+q^2=1\}.$$.

\vglue.5cm
\noindent
{\bf $\S$4. $X_n$ and monomorphisms from $\a_3$ to $\a_{n+1}$.}
\vglue.5cm
In this chapter $1\leq m\leq n$.
\vglue.5cm
\noindent
{\bf Definition.} An algebra monomorphism from $\a_m$ to $\a_n$ is a linear monomorphism $\varphi:\a_m\rightarrow\a_n$ such that 

i) $\varphi(e_0)=e_0$ (the first $e_0$ is in $\a_m$ and the second $e_0$  in $\a_n$)

ii) $\varphi(xy)=\varphi(x)\varphi(y)$ for all $x$ and $y$ in $\a_m$.

By definition we have that $\varphi(re_0)=r\varphi(e_0)$ for all $r$ in $\erre$ so $\varphi(_0\!\a_m)\subset\varphi(_0\!\a_n)$ and $\varphi(\overline{x})=\overline{\varphi(x)}$ therefore $||\varphi(x)||^2=\varphi(x)\overline{\varphi(x)}=\varphi(x)\varphi(\overline{x})=\varphi(x\overline{x})=\varphi(||x||^2)=||x||^2$ for all $x\in \a_m$ and $||\varphi(x)||=||x||$ and $\varphi$ is an orthogonal linear transformation from $\erre^{2^m-1}$ to $\erre^{2^n-1}$.

The {\it trivial} monomorphism is the one given by $\varphi(x)=
(x,0,0,\ldots,0)$ for $x\in\a_m$ and $0$ in $\a_m$

$ {\cal M}(\a_m;\a_n)$ denotes the set of algebra monomorphisms from $\a_m$ to $\a_n$.

For 
$m=n,\,\,{\cal M}(\a_m;\a_n)={\rm Aut}(\a_n)$ 
the group of algebra automorphisms of $\a_n$
\vglue.5cm
\noindent
{\bf Proposition 4.1.} ${\cal M}(\a_1;\a_n)=S(_o\!\a_n)=S^{2^n-2}$.
\vskip.3cm
\noindent
{\bf Proof:} $\a_1=\ce=Span\{e_0,e_1\}$.

If $x\in\a_1$ then $x=re_0+se_1$ and for $w\in _o\!\!\a_n$ with $||w||=1$ we have that $\varphi_w(x)=re_0+sw$ define an algebra monomorphism from $\a_1$ to $\a_n$. This can be seen by direct calculations, recalling that, Center $(\a_n)=\erre$ for all $n$ and that every associator with one real entries vanish.

Conversely, for $\varphi\in{\cal M}(\a_1;\a_n)$, set $w=\varphi(e_1)$ so $||w||=1$ and $\varphi_w=\varphi$.

\hfill Q.E.D.

\vskip.5cm
\noindent
{\bf Remark:} In particular, we have that
$$Aut(\a_1)=S^0=\ze/2=\{\hbox{\rm Identity, conjugation}\}=\{\varphi_{e_1},\varphi_{-e_1}\}.$$
To calculate ${\cal M}(\a_2;\a_n)$ for $n\geq 2$ we need to recall (see [Mo$_2$]).
\vglue.5cm
\noindent
{\bf Definition:} For $a$ and $b$ in $\a_n$. We say that 
{\it $a$ alternate with $b$}, we denote it by $a\rightsquigarrow b$, if $(a,a,b)=0$.

We say that {\it $a$ alternate strongly with $b$}, we denote it by $a\leftrightsquigarrow b$, if $(a,a,b)=0$ and $(a,b,b)=0$.

Clearly $a$ alternate strongly with $e_0$ for all $a$ in $\a_n$ and if $a$ and $b$ are linearly dependent then $a\leftrightsquigarrow b$ (by flexibility).

Also, by definition, $a$ is an alternative element if and only if $a\rightsquigarrow x$ for all $x$ in $\a_n$.

By Lemma 1.1 (1) and (2) we have that for any doubly pure element $a$ in 
$\a_n\,\, (a,a,\tilde{e}_0)=0$ and (by the above remarks) 
$\tilde{e}_0$ alternate strongly with any $a$ in $\a_n$.

For $a$ and $b$ pure elements in $\a_n,$ we define the vector subspace of $\a_n$
$$V(a;b)=Span \{e_0,a,b,ab\}.$$
{\bf Lemma 4.2.} If $(a,b)\in V_{2^n-1,2}$ and $a\leftrightsquigarrow b$ then $V(a;b)=\a_2=\ache$ the quaternions.

\vglue.5cm
\noindent
{\bf Proof:} Suppose that $(a,b)\in V_{2^n-1,2}$ and that $(a,a,b)=0$ then we have
\begin{eqnarray*}
\langle ab,a\rangle&=&\langle b,\overline{a}a\rangle=\langle b,||a||^2e_0\rangle=||a||^2\langle b,e_0\rangle=0\\
\langle ab,a\rangle&=&\langle a,b\overline{b}\rangle=\langle a,||b||^2e_0\rangle=||b||^2\langle a,e_0\rangle=0\\
||ab||^2=\langle ab,ab\rangle&=&\langle\overline{a}(ab),
b\rangle=\langle-a(ab),b\rangle=\langle-a^2b,b\rangle\\
&=&-a^2\langle b,b\rangle=||a||^2||b||^2=1
\end{eqnarray*}
so $\{e_0,a,b,ab\}$ is an orthonormal set of vectors in $\a_n$.

Finally using also that $(a,b,b)=0$ and $ab=-ba$ we may check by direct calculation that the multiplication table of $\{e_0,a,b,ab\}$ coincides with the one of the quaternions and by the identification $e_0\mapsto e_0, a\mapsto e_1, b\mapsto 
e_2$ and $ab\mapsto e_3$ we have an algebra isomorphism between $\a_2=\ache $ and $V(a;b)$.

\hfill Q.E.D.
\vglue.5cm
\noindent
{\bf Proposition 4.3.} ${\cal M}(\a_2;\a_n)=\{(a,b)\in V_{2^n-1,2}|a\leftrightsquigarrow b\}$ for $n\geq 2$.

In particular
$$Aut(\a_2)={\cal M}(\a_2;\a_2)=V_{3,2}=SO(3)$$
and 
$${\cal M}(\a_2,\a_3)=V_{7,2}.$$
{\bf Proof.} The inclusion ``$\supset$'' follows from Lemma 4.2. Conversly suppose that $\varphi \in{\cal M}(\a_2,\a_n)$ then $\varphi(e_0)=e_0,  (\varphi(e_1),\varphi(e_2))\in V_{2^n-1,2}$ and $V(\varphi(e_1),\varphi(e_2))={\rm Im}\varphi=\ache\subset\a_n$.

Since $\a_2$ is an asociative algebra and $\a_3$ is an alternative algebra we have that $a\leftrightsquigarrow b$ for any two elements in $\a_n$ for $n=2$ or $n=3$.

\hfill Q.E.D.
\vglue.5cm
\noindent
{\bf Remark.} Recall that $\tilde{\a}_n=\{e_0,\tilde{e}_0\}^\perp=\erre^{2^n-2}$ denotes the vector subspace of doubly pure elements. Since $a\leftrightsquigarrow
\tilde{e}_0$ for any element in $\tilde{\a}_n,$ we have that, if $a\in S(\tilde{\a}_n)$ i.e., $||a||=1$ then $(a,\tilde{e}_0)\in V_{2^n-1,2}$ and the assignment $a\mapsto (a,\tilde{e}_0)$ defines an inclusion from 
$S(\tilde{\a}_n)=S^{2^n-3}\hookrightarrow{\cal M}(\a_2;\a_n)\subset V_{2^n-1,2}$ which resembles ``the bottom cell'' inclusion in $V_{2^n-1,2}$.
\vglue.5cm

Now we show that $X_n$ can be identified with a subset of ${\cal M}(\a_3;\a_{n+1})$ for $n\geq 4$.
\vglue.5cm
\noindent
{\bf Lemma 4.4.} For $\alpha\in\ache^\perp_\varepsilon\subset\a_{n+1}$ and $n\geq 4$.

\noindent
(1) If $||\alpha||=1$ then $\o_\alpha={\rm Span}\{e_0,\widetilde{\varepsilon},\varepsilon,\widetilde{e}_0,\widetilde{\alpha},\alpha\varepsilon,\widetilde{\varepsilon}\alpha,\alpha\}$ is isomorphic as algebra to $\a_3$ if and only if $(\alpha,\alpha,\varepsilon)=0$.

\noindent
(2) If $\alpha=(a,b)\in\widetilde{\a}_n\times\widetilde{\a}_n$ then $(\alpha,\alpha,\varepsilon)=(0,-(a,\widetilde{e}_0,b))\in\widetilde{\a}_n\times\widetilde{\a}_n$
\vglue.5cm
\noindent
{\bf Proof.} (1) By definition $\widetilde{e}_0=e_{2^n}$ and $\varepsilon=e_{2^{n-1}}$ are elements in the canonical basis so they are alternative elements (See [Sch]). Since $\ache_\alpha$ is associative for all $\alpha$ then $(\alpha,\alpha,\widetilde{e}_0)=0$.

Clearly if $\o_\alpha\cong \a_3$ then $(\alpha,\alpha,\varepsilon)=$ because $\o_\alpha$ is an alternative algebra. 

Conversely, assume that $(\alpha,\alpha,\varepsilon)=0.$ 

We have the following multiplication table that under the mapping 

$e_0\mapsto e_0; e_1\mapsto\widetilde{\varepsilon};e_2\mapsto\varepsilon;e_3\mapsto\widetilde{e}_0;e_4\mapsto\widetilde{\alpha}; e_5\mapsto\alpha\varepsilon; e_0\mapsto\widetilde{\varepsilon}\alpha$ and $e_7\mapsto\alpha$

becomes an algebra monomorphism from $\a_3$ into $\a_{n+1}$
$$
\begin{tabular}{r|r|r|r|r|r|r|r}
$e_0$&$\widetilde{\varepsilon}$&$\varepsilon$&$\widetilde{e}_0$&
$\widetilde{\alpha}$&$\alpha\varepsilon$&$\widetilde{\varepsilon}
\alpha$&$\alpha$\\ \hline
$\widetilde{\varepsilon}$&$-e_0$&$\widetilde{e}_0$&$-\varepsilon$&$\alpha\varepsilon$&$-
\widetilde{\alpha}$&$-\alpha$&$\widetilde{\varepsilon}\alpha$\\
$\varepsilon$&$-\widetilde{e}_0$&$-e_0$&$\widetilde{\varepsilon}$&$\widetilde{\varepsilon}\alpha$&$\alpha$&$-\widetilde{\alpha}$&$-\alpha\varepsilon$\\
$\widetilde{e}_0$&$\varepsilon$&$-\widetilde{\varepsilon}$&$-e_0$&$\alpha$&$-\widetilde{\varepsilon}\alpha$&$+\alpha\varepsilon$&$-\widetilde{\alpha}$\\
$\widetilde{\alpha}$&$-\alpha\varepsilon$&$-\widetilde{\varepsilon}\alpha$&$-\alpha$&$-e_0$&$\widetilde{\varepsilon}$&$\varepsilon$&$-\widetilde{e}_0$\\
$\alpha\varepsilon$&$\widetilde{\alpha}$&$-\alpha$&$\widetilde{\varepsilon}\alpha$&$-\widetilde{\varepsilon}$&$-e_0$&$-\widetilde{e}_0$&$\varepsilon$\\
$\widetilde{\varepsilon}\alpha$&$\alpha$&$\widetilde{\alpha}$&$-\alpha\varepsilon$&$-\varepsilon$&$\widetilde{e}_0$&$-e_0$&$-\widetilde{\varepsilon}$\\
$\alpha$&$-\widetilde{\varepsilon}\alpha$&$\alpha\varepsilon$&
$\widetilde{\alpha}$&$\widetilde{e}_0$&$-\varepsilon$&$\widetilde{\varepsilon}$&$-e_0$
\end{tabular}
$$
Notice that this table is skew-symmetric with $-e_0$'s along the diagonal.

The nontrivial calculations are:

$\widetilde{\varepsilon}\widetilde{\alpha}=\alpha\varepsilon$ (by Lemma 1.1 (5)).

$\widetilde{\varepsilon}(\alpha\varepsilon)=
-\widetilde{\varepsilon\alpha)}=\widetilde{\varepsilon(\varepsilon\alpha)}=\widetilde{\varepsilon^2\alpha}=-\widetilde{e_0\alpha}=-\widetilde{\alpha}$ (by Lemma 1.1 (3)).

$\varepsilon(\widetilde{\varepsilon}\alpha)=-\varepsilon(\widetilde{\varepsilon\alpha)}=-\widetilde{\varepsilon}(\varepsilon\alpha)=
\widetilde{\varepsilon(\varepsilon\alpha)}=
-\widetilde{\alpha}$ (by Lemma 1.1 (1) and (6).

$\widetilde{\alpha}(\alpha\varepsilon)=-\widetilde{\alpha(\alpha\varepsilon)}=-\widetilde{\alpha^2\varepsilon}=\widetilde{\varepsilon}$ because $(\alpha,\alpha,\varepsilon)=0$ and $||\alpha||=1$.

$\widetilde{\alpha}(\widetilde{\varepsilon}\alpha)=-\widetilde{\alpha}(\widetilde{\varepsilon\alpha)}=+\widetilde{\widetilde{\alpha}}(\alpha\varepsilon)=-\alpha(\alpha\varepsilon)=\varepsilon$,

so we are done with (1).

To prove (2) we perform similar calculation as in Theorem 3.4.
\begin{eqnarray*}
\alpha(\alpha\varepsilon)&=&(a,b)[(a,b)(\widetilde{e}_0,0)]=(a,b)(\widetilde{a},-\widetilde{b})=(a\widetilde{a}-\widetilde{b}b,-\widetilde{b}a-b\widetilde{a})\\
&=&(-||a||^2\widetilde{e}_0-||b||^2e_0,-(b\widetilde{e}_0)a+b(\widetilde{e}_0a))\\
&=&-||\alpha||^2\varepsilon-(0,(b,\widetilde{e}_0,a))
\end{eqnarray*}
Therefore
\begin{eqnarray*}
(\alpha,\alpha,\varepsilon)&=&\alpha^2\varepsilon-\alpha(\alpha\varepsilon)=-||\alpha||^2-\alpha(\alpha\varepsilon)\\
&=&(0,(b,\widetilde{e}_0,a))=-(0,(a,\widetilde{e}_0,b))\quad{\rm by\,\,flexibility.}
\end{eqnarray*}

\hfill { Q.E.D.}
\vglue.5cm
\noindent
{\bf Notation:} For $n\geq 4$ consider the following subsets of $\a_{n+1}$
\begin{eqnarray*}
\e_n&=&\{\alpha\in\ache^\perp_\varepsilon|(\alpha,\alpha,\varepsilon)=0\},\\
S(\e_n)&=&\{\alpha\in\e_n : ||\alpha||=1\},\\
P(n)&=&\{(a,b)\in\widetilde{\a}_n\times\widetilde{\a}_n|a\,{\rm and}\,b\,{\rm are}\, \ce- {\rm collinear}\},\\
\overline{X}_n&=&\{(a,b)\in\widetilde{\a}_n\times\widetilde{\a}_n|a\neq 0,b\neq 0\,{\rm and}\, ab=0\,\}
\end{eqnarray*}
and also the following subset of monomorphisms
$${\cal M}_2(\a_3;\a_{n+1})=\{\varphi\in{\cal M}(\a_3;\a_{n+1})|\ache_\varepsilon\subset {\rm Im}\varphi\}.$$
{\bf Remark:} By Lemma 4.4 (1)  we may identify $S(\e_n)$ and ${\cal M}_2(\a_3;\a_{n+1}),$ that is, there is a one to one correspondence between this two sets.
\vglue.5cm
\noindent
{\bf Theorem 4.5.} For $n\geq 4$
\begin{itemize}
\item[(i)] $P(n)$ and $\overline{X}_n$ are subsets of $\e_n$ with $P(n)\cap\overline{X}_n=\Phi.$
\item[(ii)] There is a continous retract
$$R:\e_n\backslash P(n)\rightarrow\overline{X}_n.$$
\end{itemize}
\vglue.5cm
\noindent
{\bf Proof:} If $(a,b)\in P(n)$ then $b\in\ache_a$ or $a\in\ache_b$ and $(a,\widetilde{e}_0,b)=0$ (recall that $\ache_a$ and $\ache_b$ are associative), so by Lemma 4.4 (2)
$(\alpha,\alpha,\varepsilon)=0$ in $\a_{n+1}$ when $\alpha=(a,b)$ 
so $P(n)\subset \e_n$.

On the other hand if $(a,b)\in\overline{X}_n$ then $ab=0$ and $b\in\ache^\perp_a\subset\a_n$ and by Lemma 1.1 (1), (6) and (3)
$$(a,\widetilde{e}_0,b)=(a\widetilde{e}_0)b-a(\widetilde{e}_0b)=\widetilde{a}b+a\widetilde{b}=\widetilde{a}b+\widetilde{a}b=2\widetilde{a}b=-2\widetilde{ab}.$$
Therefore if $ab=0$ then $(a,\widetilde{e}_0,b)=0$ and $(\alpha,\alpha,\varepsilon)=0$ in $\a_{n+1}$ for $\alpha=(a,b)$ by Lemma 4.4 (2) so $\overline{X}_n\subset\e_n$.

Now if $(a,b)\in P(n)\cap \overline{X}_n$ then $b\in\ache_a$ and $ab=0$, but $\ache_a$ is associative and $a=0$ or $b=0$ which is a contradiction, therefore $P(n)\cap \overline{X}_n$ is the empty set
and we are done with (i).

To prove (ii) suppose that $\alpha=(a,b)\in\widetilde{\a}_n\times \widetilde{\a}_n$ with $a\neq 0$. Since 
$$\a_n=\ache_a\oplus\ache^\perp_a$$ 
then there are unique elements $c$ and $d$ in $\ache_a$ and $\ache^\perp_a$ respectively such that $b=c+d$. Now
$$(a,\widetilde{e}_0,b)=(a,\widetilde{e}_0,c+d)=(a,\widetilde{e}_0,c)+(a,\widetilde{e}_0,d)=0+(a,\widetilde{e}_0,d)$$
because $\ache_a$ is associative.

But by Lemma 1.1 (1), (6) and (3)
$$(a,\widetilde{e}_0,d)=(a\widetilde{e}_0)d-a(\widetilde{e}_0d)=\widetilde{a}d+a\widetilde{d}=\widetilde{ad}+\widetilde{a}d=-2\widetilde{ad}.$$
Therefore $(a,\widetilde{e}_0,b)=-2\widetilde{ad}$.

Suppose that $\alpha=(a,b)$ is in $\e_n\backslash P(n)$ then $a\neq 0$,$b\neq 0$, $$b=c+d\in\ache_a\oplus\ache^\perp_a$$
 with $d\neq 0$ and $(a,\widetilde{e}_0,b)=0$ by Lemma 4.6 (2).Thus we have that $ad=0$.

Let us define $R:\e_n\backslash P(n)\rightarrow\overline{X}_n$ as 
$R(a,b)=(a,d)$.

Then $R(a,b)=(a,b)$ if $(a,b)\in\overline{X}_n$ and $R$ is continuous, because it is the restriction of the projection map
$$
\begin{array}{ccccc}
\widetilde{\a}_n\times\a_n&\rightarrow&\widetilde{\a}_n\times(\ache_a\oplus\ache^\perp_a)&\rightarrow&\widetilde{\a}_n\times\ache^\perp_a\\
(a,b)&\rightarrow&(a,c+d)&\rightarrow&(a,d)
\end{array}
$$
which is obviously continuous.

\hfill{\bf Q.E.D.}
\vglue.5cm
\noindent
{\bf Remarks:} (1) Recall that $\widetilde{\a}_n$ is a complex vector space by making $ia=\widetilde{a}$.

By definition $(a,b)\in P(n)$ if and only if $a$ and $b$ are $\ce$-collinear for $(a,b)\in\widetilde{\a}_n\times\widetilde{\a}_n$. So
$$P(n)\cong (((\widetilde{\a}_n\backslash \{0\})\times \ce )\cup \widetilde{\a}_n).$$
(2) Consider the map $w_n:\widetilde{\a}_n\times\widetilde{\a}_n\rightarrow\widetilde{\a}_n$ given by
$$w_n(a,b)=(a,\widetilde{e}_0,b).$$
Since every associator is a pure element  $(a,\widetilde{e}_0,b)\perp\widetilde{e}_0$ because $(a,-,b)$ is a skew-symmetric linear transformation (see [Mo$_1$]) then $(a,\widetilde{e}_0,b)\in\widetilde{\a}_n$.

Now $w_n$ is a polynomia map, actually is a quadraic map, and $\e_n=w^{-1}_n(0)$ is a real algebraic set
 in $\widetilde{\a}_n\times\widetilde{\a}_n=\ache^\perp_\varepsilon=\erre^{2^{n+1}-4}$ with $0$
 in $\widetilde{\a}_n$ a \underline{singular}  value, by Lemma 4.4 (2).
 
 (3) $X_n$ is a contraction of $\overline{X_n}$ via normalization on each coordinate.

%%%%%%%%%%%%%%%%%%%%%%%%%%%%%%%%%%%%%%

\newpage
\noindent
{\bf References}

\begin{itemize}
\item[{\rm [A]}] J. Adem. Construction of some normed maps. Bol. Soc. Mat. Mexicana, (2){\bf 20} 1975, 59-79.
\item[{\rm [C$_1$]}] F. Cohen. A course in some aspects of classical homotopy theory. Lecture Notes in Mathematics 1139, 1989, Springer-Verlag, p. 11-120.
\item[{\rm [C$_2$]}]F. Cohen.  On Whitehead Squares, Cayley-Dickson algebras and rational functions. Bol. Soc.  Mat. Mexicana. (2){\bf 37} 1992, 55-62.
\item[{\rm [D]}]L.E. Dickson. On quaternions and their generalization and the history of the 8 square theorem. Annals of Math; {\bf 20} 155-171, 1919.
\item[{\rm [E-K]}]Eakin-Sathaye. On automorphisms and derivations of Cayley-Dickson algebras. J. Pure Appl. Algebra {\bf 129}, 263-280, 1990.
\item[{\rm [L]}] K. Lam. Construction of Non-singular bilinear maps. Topology 6. 1967, 423-426.
\item[{\rm [H-Y]}]S.H. Khalil and P. Yiu. The Cayley-Dickson algebras: A theorem of Hurwitz and quaternions. Bol. Soc. Sci. Lett. L\'odz Vol. XLVIII 117-169. 1997.
\item[{\rm [Mo$_1$]}]G. Moreno. The zero divisors of the Cayley-Dickson algebras over the real numbers. Bol. Soc. Mat. Mexicana (2){\bf 4}  13-27. 1998
\item[{\rm [Mo$_2$]}]G. Moreno. Alternative elements in the Cayley-Dickson algebras.\\
 hopf.math.purdue.edu/pub/Moreno.2001
\item[{\rm [Sch]}] R.D. Schafer. On the algebras formed by the Cayley-Dickson process. Amer. J. Math. 76 (1954) 435-445.
\item[{\rm [Wh]}] G. Whitehead. Elements of Homotopy theory. Graduate Texts in Math. 61 Springer-Verlag.
\end{itemize}

Departamento de Matem\'aticas.

CINVESTAV DEL IPN.

Apartado Postal 14-740.

M\'exico D.F. 

MEXICO.

gmoreno@math.cinvestav.mx

\end{document}